\newcommand{\R}{I\!\!R}

\newcommand{\Catsimp}{\underline{\mathfrak C}}
\newcommand{\Catabst}{\underline{\mathfrak D}}
\newcommand{\funtor}{\underline{\mathcal F}}

\newcommand{\dd}{\mathrm{d}}

\newcommand{\Spl}{\mathrm{Sp}}
\newcommand{\spl}{\mathrm{sp}}
\newcommand{\GL}{\mathrm{GL}}
\newcommand{\gl}{\mathrm{gl}}
\newcommand{\Id}{\mathrm{Id}}
\newcommand{\Lin}{\mathrm{Lin}}
\newcommand{\Bil}{\mathrm{B}}
\newcommand{\Bsym}{\mathrm{B}_{\textrm{sym}}}

\newcommand{\Dim}{\mathrm{dim}}

\documentclass[oneside,letterpaper,draft,10pt]{amsart}

\usepackage{amsmath}               
\usepackage{amsthm}                
\usepackage{times}
\usepackage[all]{xy}
\usepackage{amscd}

\numberwithin{equation}{section}

\title[{\tiny Conjugate Points along semi-Riemannian Geodesics}]%
{On the Distribution of Conjugate Points along semi-Riemannian Geodesics}
\author[P.\ Piccione]{Paolo Piccione}
\address{Departamento de Matem\'atica,\hfill\break\indent  Universidade de S\~ao Paulo, Brazil}
\email{piccione@ime.usp.br}
\urladdr{http://www.ime.usp.br/\~{}piccione}
\author[D.\ Tausk]{Daniel V.\ Tausk}
\address{Departamento de Matem\'atica,\hfill\break\indent  Universidade de S\~ao Paulo, Brazil}
\email{tausk@ime.usp.br}
\urladdr{http://www.ime.usp.br/\~{}tausk}
\thanks{The first author is partially sponsored by CNPq (Processo n.\ 301410/95).}

\date{October 2000}

\begin{document}


\theoremstyle{plain}\newtheorem{teo}{Theorem}[section]
\theoremstyle{plain}\newtheorem{prop}[teo]{Proposition}
\theoremstyle{plain}\newtheorem{lem}[teo]{Lemma} 
\theoremstyle{plain}\newtheorem{cor}[teo]{Corollary} 
\theoremstyle{definition}\newtheorem{defin}[teo]{Definition} 
\theoremstyle{remark}\newtheorem{rem}[teo]{Remark} 
\theoremstyle{plain} \newtheorem{assum}[teo]{Assumption}
\theoremstyle{definition}\newtheorem{example}[teo]{Example}


\begin{abstract}
Helfer in \cite{Hel} was the first to produce an example of 
a spacelike Lorentzian geodesic with a continuum of conjugate points.
In this paper we show the following result: given an interval
$[a,b]$ of $\R$ and any closed subset $F$ of $\R$ contained in $\left]a,b\right]$,
then there exists a Lorentzian manifold $(M,g)$ and a spacelike geodesic $\gamma:[a,b]\to M$ 
such that $\gamma(t)$ is conjugate to $\gamma(a)$ along $\gamma$ iff $t\in F$.
\end{abstract}

\maketitle

\begin{section}{Introduction}
\label{sec:intro}
It is well known that, in Riemannian geometry, the set of conjugate (or, more generally, focal)
points along a geodesic
is discrete; Beem and Ehrlich (see \cite{BE, BEE}) have shown that the same holds for
causal, i.e., timelike or lightlike, geodesics in a Lorentzian manifold. 
The issue of the lack of discreteness for the set of conjugate points along
a geodesic in a semi-Riemannian manifold with metric of arbitrary index has
been somewhat ignored or overlooked in the literature (see for instance
\cite[Exercise~8, pag.\ 299]{ON}, or \cite[The Index Theorem]{MasielloBook}). 
However, without a suitable nondegeneracy assumption,  the classical proof 
of discreteness for the Riemannian case does not work in the general case, and
Helfer in \cite{Hel} gave the first counterexample to the discreteness of conjugate points
along a spacelike Lorentzian geodesic. In \cite[Section~11]{Hel} it is produced an example
of a whole segment of conjugate points.

The occurrence of an infinite number of conjugate points along a compact segment of 
a semi-Riemannian geodesic is a rather pathological phenomenon, for instance it cannot happen
if the metric is real-analytic; moreover, the nondegeneracy assumption mentioned above
is {\em generic\/} (see for instance \cite{PT2}).
Nevertheless, in order to fully understand the theory of conjugate points for
non positive definite metrics, it is a natural question to ask what are the possible
``shapes'' for the set of conjugate points along a geodesic. In this paper we answer
this question by reducing the problem to the study of intersection theory of curves
in the Lagrangian Grassmannian of a symplectic space.

Given a  geodesic $\gamma:[a,b]\to M$ in a semi-Riemannian manifold $(M,g)$, 
the set $l(t)$ of pairs $\big(J(t),gJ'(t)\big)$, where $J$ is a Jacobi field 
along $\gamma$ with $J(a)=0$, is a Lagrangian subspace of the symplectic space 
$T_{\gamma(t)}M\oplus T_{\gamma(t)}M^*$ endowed with its canonical symplectic 
form; the conjugate points along $\gamma$ correspond to instants $t\in\left]a,b\right]$ 
where $l(t)$ is {\em not\/} transversal to the Lagrangian subspace 
$\{0\}\oplus T_{\gamma(t)}M^*$. The use of a (parallel)  trivialization of $TM$ along
$\gamma$ allows to associate to $l$ a curve in the Lagrangian Grassmannian $\Lambda$ 
of the fixed symplectic space $\R^n\oplus{\R^n}^*\cong\R^{2n}$. Conjugate points 
along $\gamma$ correspond therefore to intersections of this curve with the subvariety 
of $\Lambda$ consisting of Lagrangians that are not transverse to $\{0\}\oplus{\R^n}^*$. 
Details of this construction can be found in \cite{Hel, MPT, PT2, catania}. The problem of 
determining precisely which
curves of Lagrangians arise from  a semi-Riemannian geodesic is a rather difficult task. A partial result in this
direction can be found in the last section of \cite{MPT}, where it is proven that a necessary condition
for a smooth curve in the Lagrangian Grassmannian $\Lambda$  to arise
from a semi-Riemannian geodesic is that it be tangent to a singular distribution
of affine planes in $\Lambda$. However, this condition alone is not sufficient, and
attempts to produce interesting examples of conjugate points along geodesics using this characterization 
lead quickly to rather involved computations.

In this paper we introduce a new procedure for constructing a curve $\xi$ in 
the Lagrangian Grassmannian $\Lambda$
starting from a semi-Riemannian geodesic $\gamma$. This new construction is {\em canonical\/}
(see Remark~\ref{thm:remcanonical}),
i.e., it does not depend on the choice of a trivialization of $TM$ along $\gamma$,
and, again, the curve $\xi$ contains the relevant information about
the conjugate points along $\gamma$. The main feature of this new construction is that it is
very easy to characterize which curves $\xi$ actually arise from semi-Riemannian geodesics;
namely, such curves are precisely those for which $\xi'(t)$ (which is naturally identified
with a symmetric bilinear form on $\xi(t)$) is {\em nondegenerate\/} for all $t$
(Theorem~\ref{thm:ABSTRACT}).
Using this characterization, it is easy to produce examples and counterexamples
concerning the occurrence of several {\em types\/} of conjugate points along a semi-Riemannian
geodesic;  we prove in particular that any compact subset of $\R$ appears as the set of conjugate 
instants along some spacelike Lorentzian geodesic (Theorem~\ref{thm:MAIN}).

\end{section}


\begin{section}{The abstract setup}
\label{sec:setup}
Given (finite dimensional) real vector spaces $V$, $W$ we denote by $\Lin(V,W)$ the space of linear maps from $V$ to
$W$ and by $\Bil(V,W)$ the space of bilinear forms $B:V\times W\to\R$; by $\Bsym(V)$ we denote the subspace of
$\Bil(V,V)$ consisting of {\em symmetric\/} bilinear forms. The {\em index\/} of a symmetric bilinear form
$B\in\Bsym(V)$ is defined as the supremum of the dimensions of the subspaces of $V$ on which $B$ is negative
definite. We always implicitly identify the spaces
$\Bil(V,W)$ and $\Lin(V,W^*)$ by the isomorphism $B(v,w)=B(v)(w)$, where $W^*$ denotes the dual space of $W$.

Let $(M,\mathfrak g)$ be an $(n+1)$-dimensional semi-Riemannian manifold and let $\gamma:[a,b]\to M$ be a non
lightlike geodesic, i.e., $\mathfrak g(\dot\gamma,\dot\gamma)$ is not zero. Using a parallel trivialization of the
normal bundle of
$\gamma$, the Jacobi equation along $\gamma$ can be seen as a second order linear system of differential equations in
$\R^n$ of the form $v''=Rv$, where $t\mapsto R(t)$ is a smooth curve of $g$-symmetric linear endomorphisms of
$\R^n$ representing a component of the curvature tensor and $g$ is a nondegenerate symmetric bilinear form in $\R^n$
representing the semi-Riemannian metric $\mathfrak  g$ on the normal bundle of $\gamma$. An equation of the form
$v''=Rv$ with a
$g$-symmetric $R$ is called a {\em Morse--Sturm\/} system; the index of $g$ is called the
{\em index of the Morse--Sturm system}.

We recall from \cite{Hel} the following:
\begin{lem}\label{thm:lemhelfer}
Every Morse--Sturm system in $\R^n$ can be obtained by a parallel trivialization of the
normal bundle from the Jacobi equation along
a non lightlike geodesic $\gamma:[a,b]\to M$, where $(M,\mathfrak g)$ is an $(n+1)$-dimensional (conformally flat)
semi-Riemannian manifold. Moreover, the geodesic can be chosen to be either spacelike or timelike; in the first
case the index of the metric $\mathfrak  g$ equals the index of the Morse--Sturm system, and in the latter case
the index of the metric $\mathfrak  g$ equals the index of the Morse--Sturm system plus one.
\end{lem}
\begin{proof}
Consider $M=\R^{n+1}$ with coordinates $(x_1,\ldots,x_{n+1})$ and let $\gamma:[a,b]\to M$ be given by $
\gamma(t)=t\frac\partial{\partial x_{n+1}}$; consider in $M$ the metric $\mathfrak g=e^\Omega\mathfrak g_0$,
with $\mathfrak g_0=g\pm\mathrm dx_{n+1}^2$,  and $\Omega$   given by:
\[\Omega(x_1,\ldots,x_{n+1})=\pm\sum_{i,j=1}^n {g\Big(R(x_{n+1})\frac\partial{\partial
x_i},\frac\partial{\partial x_j}\Big)}x_ix_j.\]
The choice of the sign $\pm$ in the above expressions is made according to the desired causal
character of $\gamma$. It is easily checked that the Christoffel symbols of the Levi--Civita
connection of $\mathfrak g$ in the canonical basis vanish along $\gamma$; this implies that
$\gamma$ is a geodesic and that $\big(\frac\partial{\partial x_i}\big)_{i=1}^n$ gives a parallel
trivialization of the normal bundle $\dot\gamma^\perp$.
\end{proof}
 Setting $\alpha=gv'$ the Morse--Sturm equation $v''=Rv$ is written as the
following first order linear systems of differential equations:
\begin{equation}\label{eq:firstorder}
\left\{
\begin{aligned}
v'&=g^{-1}\alpha,\\
\alpha'&=gRv.
\end{aligned}\right.
\end{equation}
The coefficient matrix $\begin{pmatrix}0&g^{-1}\\gR&0\end{pmatrix}$ of \eqref{eq:firstorder} is easily seen to be a
curve in the Lie algebra $\spl(2n,\R)$ of the symplectic group $\Spl(2n,\R)$ of $\R^n\oplus{\R^n}^*$
endowed with the canonical symplectic form:
\begin{equation}\label{eq:omegacan}
\omega\big((v_1,\alpha_1),(v_2,\alpha_2)\big)=\alpha_2(v_1)-\alpha_1(v_2).
\end{equation}
Recall indeed that the Lie algebra $\spl(2n,\R)$ consists of all the matrices of the form:
\begin{equation}\label{eq:XABC}
X=\begin{pmatrix}A&B\\C&-A^*\end{pmatrix},
\end{equation}
where $A\in\Lin(\R^n)$, $B\in\Bsym({\R^n}^*)$ and $C\in\Bsym(\R^n)$.
The considerations above motivate the following:
\begin{defin}\label{thm:defsimplsist}
Let $X:[a,b]\to\spl(2n,\R)$ be a smooth curve in $\spl(2n,\R)$ and denote by $A,B,C$ the $n\times n$ blocks of $X$
as in \eqref{eq:XABC}. The system
\begin{equation}\label{eq:sistdif}
\left\{\begin{aligned}
v'&=Av+B\alpha,\\
\alpha'&=Cv-A^*\alpha,
\end{aligned}\right.
\end{equation}
is called a {\em symplectic differential system\/} in $\R^n$. With little abuse of terminology we identify the
coefficient matrix $X$ with the system \eqref{eq:sistdif} and call $X$ a symplectic differential system in $\R^n$. We
call the system $X$ {\em nondegenerate\/} if the matrix $B(t)$ is invertible for every $t\in[a,b]$; in this case, the
{\em index\/} of $X$ is defined as the index of $B(t)$ (which does not depend on $t$).

An instant $t\in\left]a,b\right]$ is said to be {\em conjugate\/} for $X$ if there exists a non zero solution
$(v,\alpha)$ of $X$ with $v(a)=v(t)=0$.
\end{defin}

The {\em fundamental matrix\/} of $X$ is the curve $[a,b]\ni t\mapsto\Phi(t)$ in the general linear group of $\R^n$
characterized by the matrix differential equation
\begin{equation}\label{eq:defPhi}
\Phi'=X\Phi,
\end{equation}
with initial condition  $\Phi(a)=\Id$; if $(v,\alpha)$ is a
solution of $X$ we have $\Phi(t)\big(v(a),\alpha(a)\big)=\big(v(t),\alpha(t)\big)$ for all $t\in[a,b]$. The fact that
$X$ takes values in $\spl(2n,\R)$ implies that $\Phi$ is actually a curve in the symplectic group $\Spl(2n,\R)$. We
will denote by $L_0$ the subspace:
\[L_0=\{0\}\oplus{\R^n}^*\subset\R^n\oplus{\R^n}^*;\]
clearly, $t\in\left]a,b\right]$ is conjugate for $X$ iff $\ell(t)\cap L_0\ne\{0\}$ where $\ell(t)$ is the subspace:
\begin{equation}\label{eq:ellt}
\ell(t)=\Phi(t)(L_0)\subset\R^n\oplus{\R^n}^*.
\end{equation}

We now define the following notion of isomorphism for symplectic differential systems.
\begin{defin}\label{thm:defisoXtildeX}
Let $X$ and $\tilde X$ be symplectic differential systems in $\R^n$. An {\em isomorphism\/} from $X$ to $\tilde X$ is
a smooth curve $\phi:[a,b]\to\Spl(2n,\R)$ with $\phi(t)(L_0)=L_0$ for all $t\in[a,b]$ satisfying either one of the
following equivalent conditions:
\begin{enumerate}
\item\label{itm:iso1} $\tilde\Phi(t)\phi(a)=\phi(t)\Phi(t)$ for all $t\in[a,b]$, where $\Phi$ and $\tilde\Phi$ denote
respectively the fundamental matrices of $X$ and $\tilde X$;

\item\label{itm:iso2} $\tilde X(t)=\phi(t)X(t)\phi(t)^{-1}+\phi'(t)\phi(t)^{-1}$ for all $t\in[a,b]$.

\end{enumerate}
If $\phi$ is an isomorphism from $X$ to $\tilde X$ we write $\phi:X\cong\tilde X$ and we say that $X$ and $\tilde X$
are {\em isomorphic\/}.
\end{defin}
It follows easily from condition \eqref{itm:iso1} above that isomorphic symplectic systems have the same conjugate
instants. Observe that an isomorphism $\phi:X\cong\tilde X$ can be written in block matrix notation as:
\[\phi=\begin{pmatrix}Z&0\\{Z^*}^{-1}W&{Z^*}^{-1}\end{pmatrix},\]
with $Z(t)\in\Lin(\R^n)$ invertible and $W(t)\in\Bsym(\R^n)$ symmetric for all $t\in[a,b]$. A straightforward
computation shows that condition
\eqref{itm:iso2} above is equivalent to:
\begin{align}
\label{eq:tildeA}&\tilde A=ZAZ^{-1}-ZBWZ^{-1}+Z'Z^{-1},\\
\label{eq:tildeB}&\tilde B=ZBZ^*,\\
&\tilde C={Z^*}^{-1}(WA+C-WBW+A^*W+W')Z^{-1},
\end{align}
where ${}^*$ denotes transposition. It follows immediately that, if $X$ is isomorphic to $\tilde X$ then $X$ is
nondegenerate iff $\tilde X$ is nondegenerate and that the indexes of $X$ and $\tilde X$ coincide.

Observe that we have a {\em category\/} $\Catsimp$ whose objects are symplectic differential systems and whose set of
morphisms from $X$ to $\tilde X$ are the isomorphisms $\phi:X\cong\tilde X$; composition of morphisms is defined in
the obvious way. Observe also that in this category {\em every morphism is an isomorphism\/}.

The study of symplectic differential systems has an interest on its own, due to the fact that such systems are
naturally in connection with solutions of Hamiltonian systems in symplectic manifolds
(see \cite{PT2}); the notion
of symplectic differential system also appears in the theory of mechanical systems subject to non holonomic
constraints and in sub-Riemannian geometry (see Section~\ref{sec:final}). In this article we are
interested in the subcategory of $\Catsimp$ consisting of {\em Morse--Sturm systems}; we say that a nondegenerate
symplectic differential system $X$ with $n\times n$ blocks $A,B,C$ is a Morse--Sturm system if $B$ is constant and
$A=0$. As we have observed in the beginning of the section, such systems always arise from the Jacobi equation along
a non lightlike semi-Riemannian geodesic by a parallel trivialization of the normal bundle. In the following lemma we
show that the category of symplectic differential systems is not ``essentially larger'' than the subcategory of
Morse--Sturm systems:
\begin{lem}\label{thm:notlarger}
Every nondegenerate symplectic differential system $X$ is isomorphic to a Morse--Sturm system.
\end{lem}
\begin{proof}
It follows easily from \eqref{eq:tildeB} that every nondegenerate symplectic differential system is isomorphic to one
whose component $B$ is constant. We
may thus assume without loss of generality that
$B$ is constant (and nondegenerate). To conclude the proof we must exhibit a smooth curve $Z$ in the Lie group
\[G=\big\{Z\in\GL(n,\R):ZBZ^*=B\big\}\]
and a smooth curve $W$ of symmetric $n\times n$ matrices such that the righthand side of \eqref{eq:tildeA} vanishes.
It suffices to take $W=\frac12\big(B^{-1}A+A^*B^{-1}\big)$ and $Z$ to be the solution of $Z'=Z(BW-A)$ with
$Z(a)=\Id$. In order to see that $Z$ takes values in $G$ simply observe that $BW-A$ is in the Lie algebra $\mathfrak
g$ of $G$ given by:
\[\mathfrak g=\big\{Y\in\gl(n,\R):YB+BY^*=0\big\}.\qedhere\]
\end{proof}

Recall that a {\em symplectic space\/} is a real finite dimensional vector space $V$ endowed with a {\em symplectic
form\/} $\omega$, i.e., $\omega$ is an antisymmetric nondegenerate bilinear form on $V$. A {\em Lagrangian subspace\/}
of $V$ is a $n$-dimensional subspace $L\subset V$ with $\omega\vert_{L\times L}=0$, where $n=\frac12\mathrm{dim}(V)$.
We denote by $\Lambda(V,\omega)$ the {\em Lagrangian Grassmannian\/} of $(V,\omega)$, i.e., the set of all Lagrangian
subspaces of $V$. The Lagrangian Grassmannian is a real-analytic compact connected $\frac12n(n+1)$-dimensional
embedded submanifold of the Grassmannian of all $n$-dimensional subspaces of $V$. We denote by $\Lambda(2n,\R)$ the
Lagrangian Grassmannian of the symplectic space $\R^n\oplus{\R^n}^*$ endowed with its canonical symplectic form
\eqref{eq:omegacan}.

Clearly, the subspace $L_0$ is Lagrangian in $\R^n\oplus{\R^n}^*$ and therefore \eqref{eq:ellt} defines a smooth
curve $\ell$ in $\Lambda(2n,\R)$; such curve is used in \cite{MPT} to study the conjugate points along a
semi-Riemannian geodesic. We now introduce the smooth curve $\xi:[a,b]\to\Lambda(2n,\R)$ given by:
\begin{equation}\label{eq:defxi}
\xi(t)=\Phi(t)^{-1}(L_0);
\end{equation}
obviously $t\in\left]a,b\right]$ is conjugate for $X$ iff $\xi(t)$ is not transversal to $\xi(a)=L_0$. This motivates
the following:
\begin{defin}
An {\em abstract symplectic system\/} is a triple $(V,\omega,\xi)$ where $(V,\omega)$ is a symplectic space and
$\xi:[a,b]\to\Lambda(V,\omega)$ is a smooth curve in the Lagrangian Grassmannian of $(V,\omega)$. An {\em
isomorphism\/} from $(V,\omega,\xi)$ to $(\tilde V,\tilde\omega,\tilde\xi)$ is a symplectomorphism
$\sigma:(V,\omega)\to(\tilde V,\tilde\omega)$ such that $\sigma(\xi(t))=\tilde\xi(t)$ for all $t\in[a,b]$; we write
$\sigma:(V,\omega,\xi)\cong(\tilde V,\tilde\omega,\tilde\xi)$. An instant
$t\in\left]a,b\right]$ is said to be {\em conjugate\/} for $(V,\omega,\xi)$ if $\xi(t)\cap\xi(a)\ne\{0\}$.
\end{defin}
It is clear that isomorphic abstract symplectic systems have the same conjugate instants. Observe that abstract
symplectic systems and their isomorphisms form a category $\Catabst$ with composition of morphisms defined in the
obvious way; as in $\Catsimp$, all morphisms of $\Catabst$ are isomorphisms.
If $X$ is a symplectic differential system and if $\xi$ is defined in \eqref{eq:defxi} then
$\funtor(X)=(\R^n\oplus{\R^n}^*,\omega,\xi)$ is an abstract symplectic system; moreover, if $\phi:X\cong\tilde X$ is
an isomorphism then $\sigma=\phi(a)$ is an isomorphism from $\funtor(X)$ to $\funtor(\tilde X)$. The rule $\funtor$
is a {\em functor\/} from the category $\Catsimp$ to the category $\Catabst$; in addition we have the following:
\begin{lem}\label{thm:equivcat}
The functor $\funtor$ is an {\em equivalence\/} from $\Catsimp$ to $\Catabst$, i.e.:
\begin{enumerate}
\item\label{itm:eqcat1} $\funtor$ is {\em full\/} and {\em faithful\/}, i.e., given symplectic differential systems
$X$ and $\tilde X$ then $\funtor$ induces a bijection from the morphisms $\phi:X\cong\tilde X$ to the morphisms
$\sigma:\funtor(X)\cong\funtor(\tilde X)$;

\item\label{itm:eqcat2} $\funtor$ is {\em surjective on isomorphism classes\/}, i.e., given an abstract symplectic
system $(V,\omega,\xi)$ there exists a symplectic differential system $X$ such that $\funtor(X)$ is isomorphic to
$(V,\omega,\xi)$.
\end{enumerate}
\end{lem}
\begin{proof}
Part \eqref{itm:eqcat1} is obtained by straightforward verification. For part \eqref{itm:eqcat2}, we describe how to
construct the symplectic differential system $X$ from the abstract symplectic system $(V,\omega,\xi)$. Choose a
smooth curve $[a,b]\ni t\mapsto\psi(t)$ where each $\psi(t)$ is a symplectomorphism from $(V,\omega)$ to
$\R^n\oplus{\R^n}^*$ (endowed with the canonical symplectic form) such that
$\psi(t)\big(\xi(t)\big)=L_0=\{0\}\oplus{\R^n}^*$ for all
$t$. Define $X$ to be the unique symplectic differential system whose fundamental matrix $\Phi$ is given by
$\Phi(t)=\psi(t)\psi(a)^{-1}$; more explicitly, take $X(t)=\Phi'(t)\Phi(t)^{-1}$. It is easy to check that
$\sigma=\psi(a)^{-1}$ is an isomorphism from $\funtor(X)$ to $(V,\omega,\xi)$.
\end{proof}

We now want to characterize which abstract symplectic systems correspond to nondegenerate symplectic differential
systems. To this aim, we recall a couple of simple facts about the geometry of the Lagrangian Grassmannian (see
for instance \cite{Duis, MPT}). Let $(V,\omega)$ be a symplectic space. A {\em Lagrangian decomposition\/} of $V$ is
a pair $(\xi_0,\xi_1)$ of Lagrangian subspaces of $V$ such that $V=\xi_0\oplus \xi_1$; to each Lagrangian
decomposition
$(\xi_0,\xi_1)$ there corresponds a chart $\varphi_{\xi_0,\xi_1}$ defined in the open subset $\Lambda^0(\xi_1)$ of
$\Lambda(V,\omega)$ consisting of those Lagrangians that are transverse to $\xi_1$. The chart $\varphi_{\xi_0,\xi_1}$
takes values in the space $\Bsym(\xi_0)$ of symmetric bilinear forms in $\xi_0$ and is defined by:
\[\varphi_{\xi_0,\xi_1}(L)=\omega(T\cdot,\cdot)\vert_{\xi_0\times \xi_0},\quad L\in\Lambda^0(\xi_1),\]
where $T:\xi_0\to \xi_1$ is the unique linear map whose graph $\mathrm{Gr}(T)=\{v+Tv:v\in \xi_0\}$ equals $L$. The
differential $\dd\varphi_{\xi_0,\xi_1}(\xi_0)$ of the chart $\varphi_{\xi_0,\xi_1}$ at $\xi_0$ gives an isomorphism
from the tangent space $T_{\xi_0}\Lambda(V,\omega)$ to the space $\Bsym(\xi_0)$; such isomorphism does not depend
on the complementary Lagrangian $\xi_1$ to $\xi_0$ and therefore for every $L\in\Lambda(V,\omega)$ there is a {\em
natural identification\/} of the tangent space $T_L\Lambda(V,\omega)$ with the space $\Bsym(L)$.

Let $L\in\Lambda(V,\omega)$ be given and consider the {\em evaluation map\/}
$\beta_L:\Spl(V,\omega)\to\Lambda(V,\omega)$ given by $\beta_L(A)=A(L)$; using local coordinates the differential of
$\beta_L$ is easily computed as:
\begin{equation}\label{eq:difbetaL}
\dd\beta_L(A)\cdot Y=\omega(YA^{-1}\cdot,\cdot)\vert_{A(L)\times A(L)},\quad
A\in\Spl(2n,\R),\ Y\in T_A\Spl(2n,\R).
\end{equation}

Let now $X$ be a symplectic differential system and define $\xi$ as in \eqref{eq:defxi}; obviously
$\xi=\beta_{L_0}\circ\Phi^{-1}$. By \eqref{eq:difbetaL} and \eqref{eq:defPhi} we have:
\[\xi'(t)=-\omega(\Phi(t)^{-1}X(t)\Phi(t)\cdot,\cdot)\vert_{\xi(t)\times\xi(t)}=
-\omega(X(t)\Phi(t)\cdot,\Phi(t)\cdot)\vert_{\xi(t)\times\xi(t)};\]
since $\omega(X(t)\cdot,\cdot)\vert_{L_0\times L_0}=B(t)$, we see that $\xi'(t)$ is the push-forward of $-B(t)$ by
the isomorphism $\Phi(t)^{-1}:L_0\to\xi(t)$. This motivates the following:
\begin{defin}\label{thm:abstractnondeg}
An abstract symplectic system $(V,\omega,\xi)$ is called {\em nondegenerate\/} when $\xi'(t)$ is a nondegenerate
symmetric bilinear form on $\xi(t)$ for all $t$. In this case, the {\em index\/} of $(V,\omega,\xi)$ is defined as
the index of $-\xi'(t)$ (which does not depend on $t$).
\end{defin}
Clearly, nondegeneracy and indexes of abstract symplectic systems  are preserved by isomorphisms; moreover, a
symplectic differential system $X$ is nondegenerate with index $k$ iff
$\funtor(X)$ is nondegenerate with index $k$ as an abstract symplectic system.

Summarizing the results of this section, we have proven the following theorem:
\begin{teo}[abstract characterization of semi-Riemannian geodesics]\label{thm:ABSTRACT}
Let $(V,\omega,\xi)$ be a nondegenerate abstract symplectic system of index $k$, with $\mathrm{dim}(V)=2n$. Then,
there exists a $(n+1)$-dimensional semi-Riemannian manifold
$(M,g)$ and a non lightlike geodesic $\gamma:[a,b]\to M$ such that $\funtor(X)$ is isomorphic to $(V,\omega,\xi)$,
where $X$ is the Morse--Sturm system obtained from the Jacobi equation along $\gamma$ by a parallel trivialization
of the normal bundle of $\gamma$ (see \eqref{eq:firstorder}). A point $\gamma(t)$, $t\in\left]a,b\right]$ is
conjugate to $\gamma(a)$ along $\gamma$ iff $t$ is a conjugate instant for $(V,\omega,\xi)$. Moreover, $\gamma$
can be chosen to be either timelike or spacelike; the index of $g$ is equal to $k+1$ in the first case and 
to $k$ in the latter case.\qed
\end{teo}
Clearly, from a strictly technical point of view, the categorical terminology adopted
in this section is unnecessary. Nevertheless, the authors believe that the employment of
this language helps the reader in perceiving the analogies between this theory and other
situations in Mathematics\footnote{%
Here are some examples. The category of simply connected Lie groups is equivalent to the category
of real, finite-dimensional Lie algebras. The same holds for the categories of geometric simplicial complexes
and abstract simplicial complexes.} 
where categorical equivalences occur. 
\end{section}


\begin{section}{Distribution of conjugate points along a geodesic}
\label{sec:distribution}
In this section we want to construct examples of conjugate points using the characterization given in
Theorem~\ref{thm:ABSTRACT}. The idea is to construct smooth curves $\xi$ of Lagrangians of
a fixed symplectic space having everywhere nondegenerate derivative, and such that $\xi(t)$
is not transversal to a fixed Lagrangian $\xi_0$ at a prescribed set of values of the parameter
$t$. Such construction is performed using local charts $\varphi_{\xi_0,\xi_1}$ in the
Lagrangian Grassmannian; in these coordinates curves of Lagrangians are identified
with curves of symmetric bilinear forms. The main technical problem to complete
the construction is to connect smoothly $\xi$ with $\xi_0$ without violating the nondegeneracy condition
on the derivative and without creating new conjugate instants (see Proposition~\ref{thm:menoschato}).
The proof of Proposition~\ref{thm:menoschato} takes inspiration from the proof of some
elementary versions of the so-called {\em H-principle\/}   of Gromov~\cite{Gromov} by the method of
convex integration;  roughly speaking, we construct a curve satisfying a certain open differential
relation by first searching for its derivative.

We start with two technical results:

\begin{lem}\label{thm:chato1}
Let $U\subset\R^k$ be a connected open set, $u\in U$ a fixed point, $\bar\tau:[c,b]\to U$ a smooth curve and $a\in\R$,
$a<c$. Then there exists $M>0$ such that for all $\eta,\eta'>0$ there exists a smooth extension
$\tau:[a,b]\to U$ of $\bar\tau$ with the following properties:
\begin{itemize}
\item $\int_a^c\tau=u(c-a)$;
\smallskip

\item $\big\Vert\tau\vert_{[a,c]}\big\Vert_\infty=\sup_{t\in[a,c]}\Vert\tau(t)\Vert\le M$;
\smallskip

\item $\tau\vert_{[a,c-\eta]}$ is constant.
\end{itemize}
\end{lem}
\begin{proof}
Let $r>0$ be such that the open ball $B(u;r)$ of center $u$ and radius $r$ is contained in $U$ and choose a smooth
curve $\tilde\gamma:[c-1,b]\to U$ such that $\tilde\gamma(c-1)=u$ and $\tilde\gamma\vert_{[c,b]}=\bar\tau$. Set
$M=\Vert u\Vert+1+\Vert\tilde\gamma\Vert_\infty$ and choose $\varepsilon>0$ small enough such that $\varepsilon<\eta$
and
\begin{equation}\label{eq:chata}
\frac\varepsilon{c-a-\varepsilon}\Vert\tilde\gamma-u\Vert_\infty<\min\{r,1\}.
\end{equation}
Now, let $\gamma:[c-\varepsilon,b]\to U$ be a smooth non decreasing reparameterization of $\tilde\gamma$ such that
$\gamma\vert_{[c,b]}=\bar\tau$ and $\gamma\vert_{[c-\varepsilon,c-\frac\varepsilon2]}\equiv u$. Choose smooth
functions
$\phi_1,\phi_2:[a,b]\to[0,1]$ with
$\phi_1+\phi_2\equiv1$ and such that the support of $\phi_1$ is contained in $\left[a,c-\frac\varepsilon2\right[$ and
the support of
$\phi_2$ is contained in
$\left]c-\varepsilon,b\right]$. Finally set:
\[\delta=\frac{-\int_a^c\phi_2(\gamma-u)}{\int_a^c\phi_1},\]
and define $\tau=\phi_1(u+\delta)+\phi_2\gamma$. To check that such $\tau$ works observe that $\Vert\delta\Vert$ is
less than or equal to the left hand side of \eqref{eq:chata}.
\end{proof}

\begin{cor}\label{thm:corchato}
Let $\bar\sigma:[c,b]\to\Bsym(\R^n)$ be a smooth map such that $\bar\sigma(c)$ is nondegenerate, $\bar\sigma'(t)$
is nondegenerate for all $t\in[c,b]$ and such that $\bar\sigma(c)$ and $\bar\sigma'(c)$ have the same index. Then,
given $a<c$ there exists a smooth extension
$\sigma:[a,b]\to\Bsym(\R^n)$ of
$\bar\sigma$ such that
$\sigma(a)=0$, $\sigma(t)$ is nondegenerate for all $t\in\left]a,c\right]$ and $\sigma'(t)$ is nondegenerate for all
$t\in[a,b]$.
\end{cor}
\begin{proof}
Simply apply Lemma~\ref{thm:chato1} to the following objects:
\begin{itemize}
\item $U=\{B\in\Bsym(\R^n):\text{$B$ is nondegenerate and it has the same index as $\bar\sigma(c)$}\}$;
\smallskip
\item $\displaystyle u=\frac{\bar\sigma(c)}{c-a}$;
\smallskip
\item $\bar\tau=\bar\sigma'$;
\smallskip
\item $\eta>0$ is chosen small enough so that $\eta M<r$, where $r>0$ is such
that the open ball $B(\bar\sigma(c);r)$ is contained in $U$.
\end{itemize}
Finally, define $\sigma(t)=\int_a^t\tau$ for $t\in[a,b]$.
\end{proof}

\begin{prop}\label{thm:menoschato}
Let $(V,\omega)$ be a symplectic space, $\xi_0\subset V$ be a Lagrangian subspace and
$\bar\xi:[c,b]\to\Lambda(V,\omega)$ be a smooth curve such that $\bar\xi(c)\cap\xi_0=\{0\}$ and
$\bar\xi'(t)\in\Bsym(\bar\xi(t))$ is nondegenerate for all $t\in[c,b]$. Then, given
$a<c$ there exists a smooth extension $\xi:[a,b]\to\Lambda(V,\omega)$ of $\bar\xi$ such that
$\xi(a)=\xi_0$, $\xi(t)\cap\xi_0=\{0\}$ for all $t\in\left]a,c\right]$ and $\xi'(t)\in\Bsym(\xi(t))$ is nondegenerate
for all $t\in[a,b]$.
\end{prop}
\begin{proof}
Let $\xi_1$ be a Lagrangian complementary to both $\xi_0$ and $\bar\xi(c)$; it's easy to see that $\xi_1$ can be
chosen such that $\varphi_{\xi_0,\xi_1}(\bar\xi(c))$ equals any prescribed nondegenerate bilinear form on $\xi_0$. In
particular, we may assume that $\varphi_{\xi_0,\xi_1}(\bar\xi(c))$ and $\bar\xi'(c)$ have the same index. Let
$b'\in\left]c,b\right]$ be such that $\bar\xi([c,b'])$ is contained in the domain of the chart
$\varphi_{\xi_0,\xi_1}$ and define $\bar\sigma:[c,b']\to\Bsym(\xi_0)\cong\Bsym(\R^n)$ by
$\bar\sigma=\varphi_{\xi_0,\xi_1}\circ\bar\xi\vert_{[c,b']}$. The conclusion follows by an application of
Corollary~\ref{thm:corchato} to $\bar\sigma$, keeping in mind that if $\sigma=\varphi_{\xi_0,\xi_1}\circ\xi$ then:
\begin{itemize}
\item[(a)]$\xi(a)=\xi_0\Leftrightarrow\sigma(a)=0$;

\item[(b)]$\xi(t)\cap\xi_0=\{0\}\Leftrightarrow\sigma(t)\ \text{nondegenerate}$;

\item[(c)]$\xi'(t)\in\Bsym(\xi(t))$ is just a {\em push-forward\/} of $\sigma'(t)\in\Bsym(\xi_0)$ by an isomorphism
between $\xi_0$ and $\xi'(t)$.

\end{itemize}
\end{proof}
We are now ready to prove the main result of the section:
\begin{teo}\label{thm:MAIN}
Let $F\subset\left]a,b\right]$ be {\em any\/} compact subset; then there
exists a 3-dimensional Lorentzian manifold $(M,g)$ and a spacelike geodesic $\gamma:[a,b]\to M$ such that $\gamma(t)$
is conjugate to $\gamma(a)$ along $\gamma$ iff $t\in F$.
\end{teo}
\begin{proof}
By Theorem~\ref{thm:ABSTRACT}, it suffices to find an abstract symplectic system $(V,\omega,\xi)$ of index $1$ with
$\Dim(V)=4$ whose set of conjugate instants is $F$. Consider the space $V=\R^2\oplus{\R^2}^*$ endowed with the
canonical symplectic form and set $\xi_0=\{0\}\oplus{\R^2}^*$; given $c\in\left]a,\inf F\right[$, we'll construct a
smooth curve
$\bar\xi:[c,b]\to\Lambda(V,\omega)$ such that $\bar\xi'(t)$ is nondegenerate for all $t$ and
$\bar\xi(t)\cap\xi_0\ne\{0\}$ iff
$t\in F$. The desired curve $\xi:[a,b]\to\Lambda(V,\omega)$ will then be obtained by applying
Proposition~\ref{thm:menoschato}. The curve $\bar\xi$ will take values in the domain of the chart
$\varphi_{\xi_0,\xi_1}$ where $\xi_1=\R^2\oplus\{0\}$; we define $\bar\xi=\varphi_{\xi_0,\xi_1}^{-1}\circ\rho$, where
$\rho:[c,b]\to\Bsym(\xi_0)\cong\Bsym(\R^2)$ is defined\footnote{%
Identifying $\Bsym(\R^2)$ with $\R^3$, then the set of degenerate bilinear forms corresponds to a double
cone $\mathcal C$. The curve $\rho(t)$ defined above takes values in a plane $\pi$ orthogonal to the axis of the cone,
and $1-R(t)$ is the distance between $\rho(t)$ and the circle $\mathcal C\cap\pi$.}
by:
\[\rho(t)=\begin{pmatrix}1+R(t)\cos(t)&R(t)\sin(t)\\
R(t)\sin(t)&1-R(t)\cos(t)\end{pmatrix},\quad t\in[c,b],\]
and $R:[c,b]\to\left]0,+\infty\right[$ is a smooth map such that $R^{-1}(1)=F$. The condition $R(t)>0$ implies that
$\rho'(t)$ is always nondegenerate and therefore also $\xi'(t)$ is nondegenerate; moreover,
$\bar\xi(t)\cap\xi_0\ne\{0\}$ iff $R(t)=1$. The existence of the required function
$R$ follows by taking $R=1-f$ in Lemma~\ref{thm:existeR} below.
\end{proof}

\begin{lem}\label{thm:existeR}
Given a closed subset $F\subset\R$, there exists a smooth map $f:\R\to\left[0,1\right[$ such that $f^{-1}(0)=F$.
\end{lem}
\begin{proof}
Write $\R\setminus F=\bigcup_{r=1}^{+\infty}I_r$ as a disjoint union of open intervals $I_r$. For each $r\ge1$ let
$f_r:\R\to\R$ be a smooth map such that:
\begin{itemize}
\item $f_r$ is zero outside $I_r$;

\item $f_r$ is positive on $I_r$;

\item $\big\Vert f^{(i)}_r\big\Vert_\infty<2^{-r}$ for $i=0,\ldots,r$, where $f_r^{(i)}$ denotes the $i$-th derivative
of
$f_r$.

\end{itemize}
To conclude the proof set $f=\sum_{r=1}^{+\infty}f_r$.
\end{proof}

Examples of non lightlike geodesics with a prescribed set of conjugate points in higher dimensional
semi-Riemannian manifolds with metric of arbitrary index can be trivially obtained 
from  Theorem~\ref{thm:MAIN} by considering orthogonal products with a flat manifold.
On the other hand, if $\gamma$ is a spacelike geodesic in a 2-dimensional Lorentzian manifold
$(M,g)$, then $\gamma$ is a timelike geodesic in the Lorentzian manifold $(M,-g)$
with the same conjugate points. This implies that the conjugate points along a geodesic in
 a 2-dimensional semi-Riemannian manifold  are always isolated. 

\end{section}

\begin{section}{Final remarks}
\label{sec:final}
\begin{rem}\label{thm:remaparecem}
If $\gamma:[a,b]\to M$ is any geodesic (of arbitrary causal character)
in a semi-Riemannian manifold $(M,g)$, then a Morse--Sturm system can be obtained
from the Jacobi equation along $\gamma$ by a parallel trivialization of
the tangent bundle $TM$ along $\gamma$. At the beginning of Section~\ref{sec:setup}
we have  defined a Morse--Sturm system from the Jacobi equation by means
of a parallel trivialization of the {\em normal bundle\/} $\dot\gamma^\perp$ of $\gamma$.
The advantage of the latter construction is that  one has a converse to the
above construction, i.e., every Morse--Sturm system arises from the
Jacobi equation along a non lightlike semi-Riemannian geodesic (Lemma~\ref{thm:lemhelfer}).

Symplectic differential systems are more generally associated to
solutions of Hamiltonian systems in a symplectic manifold endowed with
a Lagrangian distribution (details of this construction can be found in \cite{PT2, catania}). 
To each symplectic differential system is naturally associated the notion of {\em Maslov
index}; this formalism is used in \cite{PT2} to prove a Morse
index theorem for non convex Hamiltonian systems and for semi-Riemannian geometry
(see also \cite{MPT, PT3}).  In \cite{PT2} it is also defined the notions of {\em multiplicity\/} 
and of {\em signature\/}  of a conjugate instant of a symplectic differential system;
these notions, as well as that of Maslov index,  
can be defined directly in the context of abstract symplectic
systems. In the proof of Theorem~\ref{thm:MAIN} we have constructed examples containing only
conjugate instants of multiplicity one and signature zero. However,  Theorem~\ref{thm:ABSTRACT}
and Proposition~\ref{thm:menoschato} make it an easy task to produce more
{\em exotic\/} examples of geodesics of arbitrary Maslov index and having a complicated distribution of
conjugate points  of several {\em types}.
\end{rem}
 
\begin{rem}\label{thm:remcanonical}
As mentioned in the Introduction, abstract symplectic systems are {\em canonically\/} associated
to semi-Riemannian geodesics, or more generally, to solutions of Hamiltonian systems
in a symplectic manifold endowed with a Lagrangian distribution. This is done as follows.
Let $(\mathcal M,\omega)$ be a symplectic manifold (in the geodesic case, $\mathcal M=TM^*$ is the cotangent bundle
of a semi-Riemannian manifold $(M,\mathfrak g)$),   $H$ a  possibly time dependent Hamiltonian function
on $\mathcal M$ (in the geodesic case $H(p)=\frac12\mathfrak g^{-1}(p,p)$),  $\mathfrak L\subset T\mathcal M$ a
Lagrangian distribution on $\mathcal M$ (in the geodesic case, $\mathfrak L$ is the vertical subbundle
of $TTM^*$) and $\Gamma:[a,b]\to\mathcal M$ a solution of the Hamilton equations of $H$. 
An abstract symplectic system is then obtained by considering $V=T_{\Gamma(a)}\mathcal M$ and $\xi(t)$ to be the
inverse image of $\mathfrak L_{\Gamma(t)}$  in $T_{\Gamma(a)}\mathcal M$ by the Hamiltonian flow.
\end{rem}

\begin{rem}\label{thm:remfocais}
By minor modifications of the theory presented in this paper it is also possible
to treat the case of focal points to submanifolds along an orthogonal geodesic.
To this aim, one should introduce a category of pairs $(X,\ell_0)$ where
$X$ is a symplectic differential system and $\ell_0$ is a Lagrangian subspace
of $\R^n\oplus{\R^n}^*$. The Lagrangian subspace $\ell_0$ encodes the information about
the tangent space and the second fundamental form of the initial submanifold: an instant
$t\in\left]a,b\right]$ is {\em focal for
$(X,\ell_0)$} if there exists a non zero solution $(v,\alpha)$ of $X$ with
$(v(a),\alpha(a))\in\ell_0$ and $v(t)=0$. Accordingly, abstract symplectic systems
should be replaced by quadruples  $(V,\omega,\xi,\xi_0)$, where $\xi_0$ is a Lagrangian
subspace of $(V,\omega)$. Details of this construction can be found
in \cite{PT2}.
\end{rem}

\begin{rem}\label{thm:remsub}
{\em Degenerate\/} symplectic systems (systems \eqref{eq:sistdif} with coefficient $B$ degenerate) can be used to
study  stationary points of {\em constrained\/} Lagrangian problems (see \cite{Lisboa}).
An important class of examples of these stationary points are the so-called {\em sub-Riemannian geodesics}, i.e.,
geodesics in manifolds endowed with a {\em partially defined\/} metric tensor.
Also in this case, conjugate points may accumulate along a geodesic, however, we will
show in a forthcoming paper that the set of conjugate points along a geodesic is always
a finite union of isolated points and closed intervals.
\end{rem}
\end{section}

\end{document}